# Confluent hypergeometric expansions of the confluent Heun function governed by two-term recurrence relations

T.A. Ishkhanyan[1,2], V.P. Krainov[3] and A.M. Ishkhanyan[1,2]

[1]Russian-Armenian University, Yerevan, 0051 Armenia
[2]Institute for Physical Research, NAS of Armenia, Ashtarak, 0203 Armenia
[3]Moscow Institute of Physics and Technology, Dolgoprudny, 141700 Russian Federation

We show that there exist infinitely many nontrivial choices of parameters of the single confluent Heun equation for which the three-term recurrence relations governing the expansions of the solutions in terms of the confluent hypergeometric functions $_1F_1$ and $_0F_1$ are reduced to two-term ones. In such cases the expansion coefficients are explicitly calculated in terms of the Euler gamma functions.

## 1. Introduction

The single confluent Heun equation is a remarkable equation that generalizes both the Gauss ordinary and the Kummer confluent hypergeometric equations [1-3]. Though this equation has a wide coverage in contemporary physics and mathematics (see, e.g., [1-4] and references therein; a rather large (yet, not exhaustive) list concerning mainly the general relativity and cosmology is discussed in [5]), and for this reason it has extensively been studied during the recent years (see, e.g., [6-18]), however, the theory of the equation is currently far from being satisfactory for the most of applications. A reason for this is that the power-series solutions [1-3] as well as the expansions of the solutions in terms of simpler special functions (such as the incomplete beta, incomplete gamma, Bessel, Kummer, Coulomb wave, Goursat, Appell functions, etc., see, e.g., [8-12,18-23]) are governed by three- or more-term recurrence relations for successive coefficients, so that apart from a few trivial cases the expansion coefficients are not calculated explicitly. However, we have recently shown that there exist infinitely many non-trivial particular choices of the involved parameters for which the three-term recurrence relations for power-series solutions become two-term [18]. In the present paper, we show that for the expansions of the solutions in terms of the hypergeometric functions $_1F_1$ and $_0F_1$ there also exist infinitely many particular non-trivial choices of parameters for which the recurrence relations reduce to two-term ones so that the expansion coefficients are explicitly written in terms of the Euler gamma functions.

## 2. Reductions

The confluent Heun equation written in its canonical form is given as [3]

$$\frac{d^2u}{dz^2}+\left(\frac{\gamma}{z}+\frac{\delta}{z-1}+\varepsilon\right)\frac{du}{dz}+\frac{\alpha z-q}{z(z-1)}u=0. \qquad (1)$$

**2.1** Let $\varepsilon \neq 0$. An expansion of the solution of this equation in terms of the Kummer confluent hypergeometric functions is written as [10]

$$u=\sum_{n=0}^{\infty}c_n\cdot {}_1F_1(\alpha/\varepsilon+n;\gamma+\delta+n;-\varepsilon z), \qquad (2)$$

where the expansion coefficients obey the three-term recurrence relation

$$R_n c_n + Q_{n-1} c_{n-1} + P_{n-2} c_{n-2} = 0 \qquad (3)$$

with

$$R_n = -n(\gamma+\delta+n-1), \qquad (4)$$

$$Q_n = -q+\alpha+(\gamma+\delta+\varepsilon+n-1)n, \qquad (5)$$

$$P_n = -\frac{(\delta+n)(\varepsilon n+\alpha)}{\gamma+\delta+n}. \qquad (6)$$

The expansion applies if $\varepsilon$ is not zero, $\gamma+\delta$ is not zero or a negative integer and $\alpha^2+(\gamma+\delta)^2 \neq 0$. The last condition assures that the involved confluent hypergeometric functions do not degenerate to ${}_0F_0(;-\varepsilon z) = e^{-\varepsilon z}$ for all $n$.

Our result is that this recurrence admits two-term reductions for infinitely many particular choices of the involved parameters. These reductions are achieved by the following ansatz guessed from the results of [18]:

$$c_n = \left(\frac{1}{n}\frac{\prod_{k=1}^{N+2}(a_k-1+n)}{\prod_{k=1}^{N+1}(b_k-1+n)}\right)c_{n-1}, \qquad (7)$$

where we put

$$a_1,...,a_N,a_{N+1},a_{N+2}=1+e_1,...,1+e_N,\alpha/\varepsilon,\delta, \qquad (8)$$

$$b_1,...,b_N,b_{N+1}=e_1,...,e_N,\gamma+\delta \qquad (9)$$

with parameters $e_1,...,e_N$ to be defined later ($e_1,...,e_N$ are not zero or negative integers). We note that the form of this ansatz resembles that for the coefficients of the generalized hypergeometric series ${}_{N+2}F_{N+1}(a_1,...,a_{N+2};b_1,...,b_{N+1};z)$ [24,25].

To derive the result, we note that the ratio $c_n/c_{n-1}$ is explicitly written as

$$\frac{c_n}{c_{n-1}}=\frac{(\alpha/\varepsilon-1+n)(\delta-1+n)}{(\gamma+\delta-1+n)n}\prod_{k=1}^{N}\frac{e_k+n}{e_k-1+n}. \qquad (10)$$

With this, the recurrence relation (3) is rewritten as

$$R_n \frac{(\alpha/\varepsilon-1+n)(\delta-1+n)}{(\gamma+\delta-1+n)n}\prod_{k=1}^{N}\frac{e_k+n}{e_k-1+n}+Q_{n-1}+$$
$$P_{n-2}\frac{(\gamma+\delta-2+n)(n-1)}{(\alpha/\varepsilon-2+n)(\delta-2+n)}\prod_{k=1}^{N}\frac{e_k-2+n}{e_k-1+n}=0. \quad (11)$$

Substituting now $R_n$ and $P_n$ from equations (4),(6) and cancelling the common denominator, this equation becomes

$$-(\alpha/\varepsilon-1+n)(\delta-1+n)\prod_{k=1}^{N}(e_k+n)+$$
$$+Q_{n-1}\prod_{k=1}^{N}(e_k-1+n)-\varepsilon(n-1)\prod_{k=1}^{N}(e_k-2+n)=0. \quad (12)$$

This is a polynomial equation in $n$. Since the terms proportional to $n^{N+2}$ cancel, it is of degree $N+1$. Hence, we have an equation of the form

$$\sum_{m=0}^{N+1}A_m(\gamma,\delta,\varepsilon,\alpha,q;e_1,...,e_N)n^m=0. \quad (13)$$

Then, equating to zero the coefficients $A_m$ warrants the satisfaction of the three-term recurrence relation (3) for all $n$. We thus have $N+2$ equations $A_m=0$, $m=0,1,..,N+1$, of which $N$ equations serve for determination of the parameters $e_{1,2,...,N}$ and the remaining two impose restrictions on the parameters of the confluent Heun equation (1).

One of the latter restrictions is derived by calculating the coefficient $A_{N+1}$ of the term proportional to $n^{N+1}$ which is readily shown to be $\gamma-N-1-\alpha/\varepsilon$. Hence,

$$\gamma=\alpha/\varepsilon+1+N. \quad (14)$$

The second restriction imposed on the parameters of the Heun equation is checked to be a polynomial equation of the degree $N+1$ for the accessory parameter $q$.

Resolving the recurrence (10), the coefficients of expansion (2) are explicitly written in terms of the gamma functions as ($c_0=1$)

$$c_n=\frac{\Gamma(n+\delta)\Gamma\left(n+\frac{\alpha}{\varepsilon}\right)\Gamma\left(\frac{\alpha}{\varepsilon}+\delta+1+N\right)}{n!\,\Gamma(\delta)\Gamma\left(\frac{\alpha}{\varepsilon}\right)\Gamma\left(n+\frac{\alpha}{\varepsilon}+\delta+1+N\right)}\prod_{k=1}^{N}\frac{e_k+n}{e_k},\quad n=0,1,2,.... \quad (15)$$

Here are the explicit solutions of the recurrence relation (3)-(6) for $N=0,1,2$.

$N=0$:
$$\gamma=\alpha/\varepsilon+1, \quad (16)$$
$$q=\alpha(1-\delta/\varepsilon), \quad (17)$$

$$c_n = \frac{\Gamma(n+\delta)\Gamma\left(n+\frac{\alpha}{\varepsilon}\right)\Gamma\left(\frac{\alpha}{\varepsilon}+\delta+1\right)}{n!\Gamma(\delta)\Gamma\left(\frac{\alpha}{\varepsilon}\right)\Gamma\left(n+\frac{\alpha}{\varepsilon}+\delta+1\right)}. \tag{18}$$

$N=1:$ 
$$\gamma = \alpha/\varepsilon + 2, \tag{19}$$

$$q^2 + q\left(\frac{\alpha}{\varepsilon}(1+2\delta) - 2\alpha - \varepsilon + \delta\right) + \frac{\alpha}{\varepsilon}\left(\frac{\alpha}{\varepsilon}(\delta-\varepsilon)(1+\delta-\varepsilon) + (\delta+\delta^2 - 2\delta\varepsilon + \varepsilon^2)\right) = 0, \tag{20}$$

$$e_1 = q - \alpha + \delta - \varepsilon + \frac{\alpha(1+\delta)}{\varepsilon}, \tag{21}$$

$$c_n = \frac{\Gamma(n+\delta)\Gamma\left(n+\frac{\alpha}{\varepsilon}\right)\Gamma\left(\frac{\alpha}{\varepsilon}+\delta+2\right)}{n!\Gamma(\delta)\Gamma\left(\frac{\alpha}{\varepsilon}\right)\Gamma\left(n+\frac{\alpha}{\varepsilon}+\delta+2\right)} \frac{e_1+n}{e_1}. \tag{22}$$

$N=2:$ 
$$\gamma = \alpha/\varepsilon + 3, \tag{23}$$

$$\begin{aligned}&q^3 + q^2\left(1+3\delta-3\varepsilon-3\alpha+\frac{3\alpha(1+\delta)}{\varepsilon}\right) + q\left(\frac{\alpha^2}{\varepsilon^2}\left(2+3\delta^2-6\delta(\varepsilon-1)+3\varepsilon(\varepsilon-2)\right)+\right.\\&\left.+2\frac{\alpha}{\varepsilon}\left(1+3\delta^2+\delta(5-6\varepsilon)+3\varepsilon(\varepsilon-1)\right)+2\left(\delta+\delta^2-2\delta\varepsilon+\varepsilon^2\right)\right) + \\&\frac{\alpha}{\varepsilon}\left(\frac{\alpha^2}{\varepsilon^2}(\delta-\varepsilon)(1+\delta-\varepsilon)(2+\delta-\varepsilon)+\frac{\alpha}{\varepsilon}(1+\delta-\varepsilon)\left(3\delta(2+\delta)-2(1+3\delta)\varepsilon+3\varepsilon^2\right)+\right.\\&\left.2\left(\delta(1+\delta)(2+\delta)-3\delta(1+\delta)\varepsilon+3\delta\varepsilon^2-\varepsilon^3\right)\right) = 0,\end{aligned} \tag{24}$$

$$e_1 + e_2 = q + 1 - \alpha + 2\delta - 2\varepsilon + \frac{\alpha(2+\delta)}{\varepsilon}, \tag{25}$$

$$2e_1 e_2 + \left(\varepsilon - \delta - q - \frac{\alpha}{\varepsilon}(1+\delta-\varepsilon)\right)(e_1+e_2) - \frac{\alpha+2\alpha\delta+\varepsilon(\delta+\varepsilon)}{\varepsilon} = 0, \tag{26}$$

$$c_n = \frac{\Gamma(n+\delta)\Gamma\left(n+\frac{\alpha}{\varepsilon}\right)\Gamma\left(\frac{\alpha}{\varepsilon}+\delta+3\right)}{n!\Gamma(\delta)\Gamma\left(\frac{\alpha}{\varepsilon}\right)\Gamma\left(n+\frac{\alpha}{\varepsilon}+\delta+3\right)} \frac{(e_1+n)(e_2+n)}{e_1 e_2}. \tag{27}$$

These results can readily be checked by direct verification of the recurrence (3)-(6).

Since $e_1,\ldots,e_N$ are not zero or negative integers, it follows from equation (10) that $c_n$ vanishes only if $\alpha/\varepsilon$ or $\delta$ is zero or a negative integer. The first case is not much interesting because then the resulting solution (2) is a polynomial in $z$. As regards the second choice, it can be checked that for integer $\delta < -N$ the solution identically vanishes, hence, the above recurrence generates non-zero finite-sum expansions of the solutions of the confluent Heun equation if

$$-N \leq \delta \leq 0. \tag{28}$$

It can further be checked that in these cases expansion (2) reduces to particular cases (because of the additional restriction $\gamma = \alpha/\varepsilon + 1 + N$, see equation (14)) of the $_{N+1}F_{1+N}$ solutions of the confluent Heun equation presented in [18].

**2.2** There exist expansions of the solutions of the confluent Heun equation in terms of the Kummer confluent hypergeometric functions the form of which differs from those applied in expansion (2). It can be checked, however, that not for all such expansions two-term reductions of the governing recurrences are possible. For instance, consider the following expansion (compare with (2)-(6)) [10]:

$$u = \sum_{n=0}^{\infty} c_n \cdot {}_1F_1(\alpha/\varepsilon + n; \gamma; -\varepsilon z), \tag{29}$$

where the expansion coefficients obey the three-term recurrence relation

$$R_n c_n + Q_{n-1} c_{n-1} + P_{n-2} c_{n-2} = 0 \tag{30}$$

with

$$R_n = n(\alpha/\varepsilon - \gamma + n), \tag{31}$$

$$Q_n = -q + n\gamma + (\alpha/\varepsilon + n)(\varepsilon - \delta - 2n), \tag{32}$$

$$P_n = (\alpha/\varepsilon + n)(\delta + n). \tag{33}$$

As it is seen, this expansion applies if $\varepsilon \neq 0$ and $\gamma$ is not zero or a negative integer. The possible two-term reduction of the recurrence (30) through ansatz (7) in this case reads

$$\frac{c_n}{c_{n-1}} = \frac{(\alpha/\varepsilon - 1 + n)(\delta - 1 + n)}{(\alpha/\varepsilon - \gamma + n)n} \prod_{k=1}^{N} \frac{e_k + n}{e_k - 1 + n}. \tag{34}$$

Hence, equation (30) is rewritten as

$$R_n \frac{(\alpha/\varepsilon - 1 + n)(\delta - 1 + n)}{(\alpha/\varepsilon - \gamma + n)n} \prod_{k=1}^{N} \frac{e_k + n}{e_k - 1 + n} + Q_{n-1} + P_{n-2} \frac{(\alpha/\varepsilon - \gamma - 1 + n)(n-1)}{(\alpha/\varepsilon - 2 + n)(\delta - 2 + n)} \prod_{k=1}^{N} \frac{e_k - 2 + n}{e_k - 1 + n} = 0. \tag{35}$$

Substituting $R_n$ and $P_{n-2}$ from equations (31),(33) and cancelling the common denominator, we obtain the following polynomial equation in $n$:

$$(\alpha/\varepsilon-1+n)(\delta-1+n)\prod_{k=1}^{N}(e_k+n)+ \tag{36}$$
$$+Q_{n-1}\prod_{k=1}^{N}(e_k-1+n)+(\alpha/\varepsilon-\gamma-1+n)(n-1)\prod_{k=1}^{N}(e_k-2+n)=0.$$

Since the terms proportional to $n^{N+2}$ cancel, this equation is of degree $N+1$. The point now is that the coefficient of the remaining highest-order term proportional to $n^{N+1}$ is checked to be equal to $\varepsilon$. Since $\varepsilon$ cannot be zero, this term cannot be canceled. Thus, the two-term reduction of the recurrence (30)-(33) via ansatz (7) is not possible.

**2.3** Here is another expansion of the solution of the confluent Heun equation [10] for which the two-term reductions are possible:

$$u = \sum_{n=0}^{\infty} c_n \cdot {}_1F_1(\alpha/\varepsilon;\gamma+n;-\varepsilon z), \tag{37}$$

where the expansion coefficients obey the three-term recurrence relation

$$R_n c_n + Q_{n-1} c_{n-1} + P_{n-2} c_{n-2} = 0 \tag{38}$$

with

$$R_n = -n(\gamma+\delta-1+n), \tag{39}$$

$$Q_n = -q+\alpha-\delta\varepsilon-n(1-\gamma-\delta+\varepsilon)+n^2, \tag{40}$$

$$P_n = \varepsilon\frac{(-\alpha/\varepsilon+\gamma+\delta+n)(\delta+n)}{\gamma+\delta+n}. \tag{41}$$

This expansion applies if $\varepsilon \neq 0$ and $\gamma$ is not zero or a negative integer. Acting essentially as above, we arrive at the following explicit solution of this recurrence:

$$c_n = \frac{\Gamma(n+\delta)\Gamma(n+\gamma-1-N)\Gamma(\gamma+\delta)}{n!\,\Gamma(\delta)\Gamma(\gamma-1-N)\Gamma(n+\gamma+\delta)}\prod_{k=1}^{N}\frac{e_k+n}{e_k}, \quad n=0,1,2,\ldots. \tag{42}$$

with

$$\delta = \alpha/\varepsilon - 1 - N. \tag{43}$$

The equations for the accessory parameter $q$ and auxiliary parameters $e_k$ for $N=0,1,2$ are as follows.

$N=0$:
$$q = (1-\gamma)\delta + \varepsilon, \tag{44}$$

$N=1$:
$$e_1 = q-2+\gamma-\delta+\gamma\delta-\varepsilon, \tag{45}$$

$$q^2 + q(\gamma-2-3\delta+2\gamma\delta-3\varepsilon)+ \tag{46}$$
$$(2-3\gamma+\gamma^2)\delta^2+2\varepsilon(2-\gamma+\varepsilon)+\delta(2+\gamma^2+4\varepsilon-3\gamma-3\gamma\varepsilon)=0,$$

$N=2$:
$$e_1+e_2+2q-7+3\gamma+3\delta-\varepsilon=0, \tag{47}$$

$$2e_1e_2 +(1-q-\gamma-\gamma\delta+2\delta+2\varepsilon)(e_1+e_2)+2q-7+3\gamma+3\delta-\varepsilon=0, \quad (48)$$

$$q^3 + q^2\left(3\gamma(1+\delta)-2(4+3\delta+3\varepsilon)\right)+ \\
q(12+26\delta+\gamma^2\left(2+3\delta(2+\delta)\right)+36\varepsilon+11(\delta+\varepsilon)^2- \\
2\gamma(5+7\varepsilon+\delta(13+6\delta+6\varepsilon)))+(\gamma-1)(\gamma-2)(\gamma-3)\delta(1+\delta)(2+\delta)- \\
2\big(9(1+\delta)(2+\delta)+\gamma^2\left(3+\delta(7+3\delta)\right)-\gamma(15+\delta(29+11\delta))\big)\varepsilon+ \\
\left(-18(2+\delta)+\gamma(15+11\delta)\right)\varepsilon^2-6\varepsilon^3. \quad (49)$$

**2.4** Consider now the case $\varepsilon=0$. This case should be examined separately because the nature of the irregular singularity of the confluent Heun equation at infinity is changed [26]. To start the discussion for this case, we note that if additionally $\alpha=0$ the confluent Heun equation is reduced to the Gauss ordinary hypergeometric equation. Therefore, since this specification is not interesting, we suppose $\alpha\neq 0$.

An appropriate expansion of the solution of the confluent Heun equation for $\varepsilon=0$ is constructed in terms of the hypergeometric functions $_0F_1$:

$$u = \sum_{n=0}^{\infty} c_n \cdot {_0F_1}(;\gamma+\delta+n;-\alpha z), \quad (50)$$

where the expansion coefficients obey the three-term recurrence relation

$$R_n c_n + Q_{n-1} c_{n-1} + P_{n-2} c_{n-2} = 0 \quad (51)$$

with
$$R_n = -n(\gamma+\delta+n-1), \quad (52)$$

$$Q_n = -q+\alpha+n(\gamma+\delta+n-1), \quad (53)$$

$$P_n = -\frac{(\delta+n)\alpha}{\gamma+\delta+n}. \quad (54)$$

The expansion is valid if $\alpha\neq 0$ and $\gamma+\delta$ is not zero or a negative integer.

We note that this recurrence is a particular case of the recurrence relation (3)-(6) achieved by putting $\varepsilon=0$. Therefore, the above technique applies for this case as well (just, we keep in mind that the expansion functions are now the hypergeometric functions $_0F_1$). The resulting two-term reduction reads ($c_0=1$)

$$\gamma = -N, \quad N=0,1,2,\ldots, \quad (55)$$

$$c_n = \frac{\alpha^n \Gamma(\delta-N)}{n!\,\Gamma(n+\delta-N)} \prod_{k=1}^{N} \frac{e_k+n}{e_k}, \quad (56)$$

For the simplest case $N=0$ this solution is valid if $q=0$. For the next two cases $N=1$ and $N=2$ the parameters should obey the following equations:

$$N=1: \qquad qe_1 + \alpha = 0, \quad q^2 + q(1-\delta) + \alpha = 0, \tag{57}$$

$$N=2: \qquad e_1 + e_2 + q + 3 - 2\delta = 0, \quad q(e_1 e_2 - \alpha) + 2\alpha(\delta-1) = 0, \tag{58}$$

$$q^3 + q^2(4-3\delta) + q(4-4\alpha-6\delta+2\delta^2) + 4\alpha(\delta-1). \tag{59}$$

These solutions are readily checked by direct verification of the recurrence (51)-(54).

## 3. Discussion

Thus, we have shown that for the expansions of the solutions of the confluent Heun equation in terms of the confluent hypergeometric functions $_0F_1$ and $_1F_1$ there exist infinitely many particular choices of equation parameters for which the three-term recurrence relations for expansion coefficients become two-term. For these cases, we have presented the explicit solutions for expansion coefficients written in terms of the Euler gamma functions. Since in all cases the computing algorithms are regular and systematic, these solutions can easily be incorporated into a computer algebra system.

The approach we applied to derive the results is essentially based on the ansatz (7). This ansatz is guessed based on the pioneering observations by Letessier [27-29] and further notable discussions by Letessier et al. [30], Maier [31] and Takemura [32]. We have already applied the ansatz to reveal several two-term reductions of the three-term recurrence relations governing the series solutions of the general and single-confluent Heun equations in terms of powers or hypergeometric functions $_0F_1$, $_1F_1$ and $_2F_1$ [27,33]. An observation here is that in all these cases the coefficients of the recurrence relations are (or can be rewritten as) polynomials, in the summation index $n$, of the degree at most 3. However, there exist other three-term expansions of the solutions of the Heun equations for which the coefficients of the recurrence relations are polynomials of higher degree (e.g., the expansions by Leaver [20], Mano et al. [21], El-Jaick and Figueiredo [8], or the eirlier Gauss-hypergeometric expansions of the general Heun functions by Svartholm [34] and Erdélyi [35]). It is interesting if the ansatz (7) works for these cases too. As we have seen in section **2.2**, not all three-term recurrences admit two-term reductions. Hence, the answer is not a priori evident and these recurrences need to be examined case by case. We conclude by noting that another appealing quest is if the ansatz is capable to produce non-trivial reductions for four- and more-term recurrences. We hope to explore these possibilities in the near future.

**Acknowledgments**

The work was supported by the Armenian Science Committee (SC Grants No. 18RF-139 and 18T-1C276), the Armenian National Science and Education Fund (ANSEF Grant No. PS-4986), the Russian-Armenian (Slavonic) University at the expense of the Ministry of Education and Science of the Russian Federation, the Russian Foundation for Basic Research (project No. 18-52-05006), and the Ministry of Education and Science of Russian Federation (Grant No. 3.873.2017/4.6). Tigran Ishkhanyan acknowledges the support from the French Embassy in Armenia for a doctoral grant as well as the Agence universitaire de la Francophonie and the Armenian Science Committee for a Scientific Mobility grant.